\documentclass[a4paper,12pt]{article}

\usepackage{amsthm}
\usepackage{amsmath,amssymb,latexsym,amsfonts,mathrsfs}

\theoremstyle{plain}
\newtheorem{Thm}{Theorem}[section]

\theoremstyle{definition}
\newtheorem{Def}[Thm]{Definition}
\newtheorem{Rem}[Thm]{Remark}

\newcommand{\Proof}[2][Proof]{\begin{proof}[{#1}] #2 \end{proof}}

\renewcommand{\d}{{\rm d}} 
\renewcommand{\hat}{\widehat}
\renewcommand{\tilde}{\widetilde}

\newcommand{\dist}{\stackrel{{\rm d}}{=}}
\newcommand{\cdist}{\stackrel{{\rm d}}{\longrightarrow}}

\newcommand{\e}{{\rm e}} 
\newcommand{\tend}[2]{\mathrel{\mathop{\longrightarrow}\limits^{#1}_{#2}}}

\numberwithin{equation}{section}

\newcommand{\rbra}[1]{\left( #1 \right)} 
\newcommand{\cbra}[1]{\left\{ #1 \right\}} 
\newcommand{\sbra}[1]{\left[ #1 \right]} 

\newcommand{\bR}{\ensuremath{\mathbb{R}}}
\newcommand{\cN}{\ensuremath{\mathcal{N}}}
\newcommand{\cF}{\ensuremath{\mathcal{F}}}

\begin{document}
\begin{center}
{\Large \bf Scaling limit of d-inverse 
of Brownian motion with functional drift}
\end{center}

\begin{center}
Kouji Yano\footnote{
Graduate School of Science, Kobe University, Kobe, JAPAN.
}\footnote{The research of this author was supported by KAKENHI (20740060).} 
\quad \text{and} \quad 
Katsutoshi Yoshioka\footnotemark[1]\footnote{Current affiliation: FUJITSU Kansai Systems Ltd.}
\end{center}

\begin{center}
{\small \today}
\end{center}

\begin{abstract}
The d-inverse is a generalized notion of inverse of a stochastic process 
having a certain tendency of increasing expectations. 
Scaling limit of the d-inverse of Brownian motion with functional drift is studied. 
Except for degenerate case, the class of possible scaling limits 
is proved to consist of the d-inverses of Brownian motion without drift, 
one with explosion in finite time, and one with power drift. 
\end{abstract}


\section{Introduction}

For (general) stock price $ S=(S_t)_{t \ge 0} $, 
the European call option price 
with strike $ K $ and maturity $ t $ is given as 
\begin{align}
C(t) := E \sbra{ \max \cbra{ S_t-K , 0 } } . 
\label{}
\end{align}
Suppose that the stock price is given as the {\em geometric Brownian motion} 
with volatility $ \sigma > 0 $ and drift $ \mu \in \bR $: 
\begin{align}
\d S_t = \sigma S_t \d B_t + \mu S_t \d t 
, \quad S_0 = s_0\in (0,\infty) , 
\label{}
\end{align}
where $ B=(B_t)_{t \ge 0} $ denotes a one-dimensional standard Brownian motion. 
Letting $ \tilde{\mu} = \mu - \sigma^2/2 $, 
we have an explicit expression of $ S=S^{(\sigma,\mu)} $ as follows: 
\begin{align}
S^{(\sigma,\mu)}_t = s_0 \exp \rbra{ \sigma B_t + \tilde{\mu} t } . 
\label{eq: gbm}
\end{align}
If $ \tilde{\mu} = - \sigma^2/2 $, then we may express $ C(t) $ 
explicitly, in terms of the cumulative distribution function 
$ \cN(x) = \int_{-\infty }^x \e^{-x^2/2} \d x / \sqrt{2 \pi} $ of the standard Gaussian, 
as 
\begin{align}
C(t) = 
s_0 \cN \rbra{ - \frac{1}{\sigma \sqrt{t}} \log \frac{K}{s_0} + \frac{1}{2} \sigma \sqrt{t} } 
- 
K \cN \rbra{ - \frac{1}{\sigma \sqrt{t}} \log \frac{K}{s_0} - \frac{1}{2} \sigma \sqrt{t} } , 
\label{}
\end{align}
which is a special case of the well-known Black--Scholes formula. 
We may verify, by a direct computation, that $ C(t) $ is increasing in $ t>0 $; 
see Madan--Roynette--Yor \cite{MRY-put}. 

Note that $ S^{(\sigma,\mu)} $ is a submartingale if and only if $ \mu \ge 0 $. 
In this case, 
we can verify, without computing it explicitly, that $ C(t) $ is increasing in $ t > 0 $. 
(In this paper, we mean non-decreasing by increasing.) 
More generally, for any increasing convex function $ \varphi $, 
we may apply Jensen's inequality to see that 
\begin{align}
E[\varphi(S_s)] \le E[\varphi(E[S_t|\cF_s])] \le E[\varphi(S_t)] 
, \quad 0 < s < t . 
\label{}
\end{align}
In this sense, the submartingale property may be considered a tendency of increasing expectations. 

To characterize another tendency of increasing expectations, 
The following notion was introduced by 
Madan--Roynette--Yor \cite{MRY-opt} and was developed by 
Profeta--Roynette--Yor \cite{MR2582990}: 

\begin{Def}
Let $ R=(R_t)_{t \ge 0} $ denote a stochastic process taking values on $ [0,\infty ) $ 
defined on a measurable space equipped with a family of probability measures $ (P_x)_{x \ge 0} $. 
Suppose that $ R $ is a.s. continuous 
and such that $ P_x(R_0=x)=1 $ for all $ x \ge 0 $. 
\begin{enumerate}
\item 
$ R $ is said {\em to admit an increasing pseudo-inverse} if 
$ P_x(R_t \ge y) $ is increasing in $ t \ge 0 $ 
for all $ y>x $ 
and if $ P_x(R_t \ge y) \to 1 $ as $ t \to \infty $ for all $ y>x $. 
\item 
A family of random variables $ (Y_{x,y})_{y>x} $ 
defined on a probability space $ (\Omega,\cF,P) $ 
is called {\em pseudo-inverse} of $ R $ 
if for any $ y>x $ it holds that 
\begin{align}
P_x(R_t \ge y) = P(Y_{x,y} \le t) . 
\label{eq: pseudo-inverse definition}
\end{align}
\end{enumerate}
\end{Def}

We would like here to introduce the following alternative notion, 
which is a slight modification of the pseudo-inverse: 

\begin{Def}
Let $ x_0 \in \bR $. Let $ X=(X_t)_{t \ge 0} $ 
be a stochastic process taking values in $ [-\infty ,\infty ] $. 
\begin{enumerate}
\item 
$ X $ is called {\em d-increasing} on $ [x_0,\infty ) $ if 
$ P(X_t \ge x) $ is increasing in $ t \in (0,\infty ) $ 
for all $ x \in [x_0,\infty ) $. 
\item 
A family of random variables $ (Y_x)_{x \ge x_0} $ 
is called {\em d-inverse} of $ X $ on $ [x_0,\infty ) $
if the following assertions hold: 
\subitem { (i)} 
for any $ x \in [x_0,\infty ) $, 
the $ Y_x $ is a random variable taking values in $ [0,\infty ] $; 
\subitem {(ii)} 
for any $ x \in [x_0,\infty ) $ and for a.e. $ t \in (0,\infty ) $, it holds that 
\begin{align}
P(X_t \ge x) = P(Y_x \le t) . 
\label{eq: d-inverse definition}
\end{align}
\end{enumerate}
\end{Def}

We note that 
$ X $ is d-increasing on $ [x_0,\infty ) $ 
if and only if $ X $ admits some d-inverse $ (Y_x)_{x \ge x_0} $. 
We also note that 
if $ P(X_t \ge x) $ is right-continuous in $ t \in (0,\infty ) $, 
then the identity \eqref{eq: d-inverse definition} holds for all $ t \in (0,\infty ) $. 

If $ t \mapsto X_t $ is a.s. increasing, 
then $ X $ is d-increasing 
and its d-inverse is given by its inverse in the usual sense. 
The d-inverse may be a generalized notion of inverse 
in the sense of probability distribution. 

Let $ S $ be a stochastic process such that 
$ P(S_t \ge x) $ is right-continuous in $ t \in (0,\infty ) $. 
We note that $ S $ is d-increasing on $ [x_0,\infty ) $ if and only if 
$ E[\varphi(S_t)] $ is increasing in $ t>0 $ 
for all increasing (possibly non-convex) function $ \varphi $ 
whose support is contained in $ [x_0,\infty ) $ 
such that $ E[\varphi(S_t)] < \infty $ for all $ t>0 $. 
In fact, for the sufficiency, 
it holds that 
\begin{align}
E[\varphi(S_t)] 
= \varphi(x_0)P(S_t \ge x_0) + \int_{x_0}^{\infty } P(S_t \ge x) \d \varphi(x) 
, \quad t>0 , 
\label{}
\end{align}
which shows that $ E[\varphi(S_t)] $ is increasing in $ t>0 $; 
the necessity is obvious since 
\begin{align}
E[1_{[x,\infty )}(S_t)] = P(S_t \ge x) . 
\label{}
\end{align}
In particular, if $ S $ is a non-negative process 
such that $ P(S_t \ge x) $ is right-continuous in $ t \in (0,\infty ) $, 
then the condition that $ S $ is d-increasing on $ [0,\infty ) $ 
is stronger 
than the one that $ S $ has the same one-dimensional marginals with a submartingale; 
see Remark \ref{rem: one-submart}. 

In this paper, we confine ourselves to the class of processes of the form 
\begin{align}
B^{(\rho)}_t = B_t + \rho(t) 
\label{}
\end{align}
for some increasing function $ \rho(t) $. 
We may call $ B^{(\rho)} $ {\em Brownian motion with functional drift}. 
This process appears in {\em geometric Brownian motion with functional coefficients} 
as follows. 
Let $ \sigma(t) $ and $ \mu(t) $ be positive functions on $ [0,\infty ) $ 
and define 
\begin{align}
\d S_t = \sigma(t) S_t \d B_t + \mu(t) S_t \d t 
, \quad S_0 = s_0 > 0. 
\label{}
\end{align}
The resulting process $ S=S^{(\sigma,\mu)} $ is given in the explicit form as 
\begin{align}
S^{(\sigma,\mu)}_t = s_0 \exp \rbra{ \int_0^t \sigma(s) \d B_s + \int_0^t \tilde{\mu}(s) \d s } , 
\label{}
\end{align}
where $ \tilde{\mu}(t) = \mu(t) - \sigma^2(t)/2 $. 
If we set $ a(t) = \int_0^t \sigma(u)^2 \d u $, $ b(t) = \int_0^t \tilde{\mu}(u) \d u $ 
and set $ \rho(t) = b(a^{-1}(t)) $, 
then we obtain 
\begin{align}
S^{(\sigma,\mu)}_{a^{-1}(t)} = s_0 \exp \rbra{ \beta_t + \rho(t) } , 
\label{}
\end{align}
where $ \beta = (\beta_t)_{t \ge 0} $ denotes a new Brownian motion.

The aim of this paper is to study scaling limit of the d-inverse on $ [0,\infty ) $ 
of $ B^{(\rho)} $ for positive drift $ \rho $. 
By {\em scaling limit} of d-inverse $ Y^{(\rho)} = (Y^{(\rho)}_x)_{x \ge 0} $ of $ B^{(\rho)} $ 
we mean a process $ Z=(Z_x)_{x \ge 0} $ such that 
\begin{align}
\frac{1}{\lambda} Y^{(\phi_1(\lambda) \rho)}_{\phi_2(\lambda) x} 
\tend{\rm d}{\lambda \to 0+} Z_x 
\quad \text{for all $ x \in [0,\infty) $} 
\label{}
\end{align}
for some scaling functions $ \phi_1 $ and $ \phi_2 $. 
We assume that the ratio $ \phi_2(\lambda) / \sqrt{\lambda} $ converges to a constant 
as $ \lambda \to 0+ $. 
We shall prove that the class of possible scaling limits 
consists, except for degenerate case, 
of the d-inverses of the following processes: 
\begin{enumerate}
\item Brownian motion without drift $ B_t $; 
\item {\em Brownian motion with explosion in finite time}: 
$ B_t + \infty 1_{\{ t \ge t_0 \}} $, with $ t_0 \in (0,\infty ) $; 
\item {\em Brownian motion with power drift}: 
$ B_t + c t^{\alpha } $, with $ c \in (0,\infty) $ and $ \alpha \ge 1/2 $. 
\end{enumerate}
Cases (i) and (ii) can be obtained from (iii) by taking limits; in fact, 
Case (i) can be obtained from (ii) as $ t_0 \to \infty $ 
and Case (ii) can be obtained from (iii) 
by setting $ c=t_0^{-\alpha } $ and letting $ \alpha \to \infty $.

Here we make several remarks. 

\begin{Rem}
Monotonicity of more general option prices for more general stock processes 
have been studied 
by Hobson (\cite{MR2563207}, \cite{MR1620358}), 
Henderson--Hobson (\cite{MR1790134}, \cite{MR1978894}), 
and Kijima \cite{MR1926239}. 
\end{Rem}

\begin{Rem}
Let $ X^{(1)} $ and $ X^{(2)} $ be two random variables taking values in $ [-\infty ,\infty ] $ 
and let $ x_0 \in \bR $. 
We write 
\begin{align}
X^{(1)} \le_{\rm st} X^{(2)} 
\quad \text{on $ [x_0,\infty ) $} 
\label{}
\end{align}
if 
\begin{align}
P(X^{(1)} \ge x) \le P(X^{(2)} \ge x) 
\quad \text{for all $ x \in [x_0,\infty ) $}. 
\label{}
\end{align}
The relation $ \le_{\rm st} $ on $ [x_0,\infty ) $ 
is a partial order on the class of random variables. 
It may be called {\em usual stochastic order on $ [x_0,\infty ) $} 
(see also Shaked--Shanthikumar \cite{MR2265633}). 
We point out that a process $ (X_t)_{t \ge 0} $ is d-increasing on $ [x_0,\infty ) $ 
if and only if $ t \mapsto X_t $ is increasing in d-order on $ [x_0,\infty ) $. 
\end{Rem}

\begin{Rem} \label{rem: one-submart}
Let $ X^{(1)} $ and $ X^{(2)} $ be two random variables taking values in $ \bR $. 
We write 
\begin{align}
X^{(1)} \le_{\rm icx} X^{(2)} 
\label{}
\end{align}
if 
\begin{align}
E[\varphi(X^{(1)})] \le E[\varphi(X^{(2)})] 
\quad 
\text{for all increasing convex function $ \varphi $}. 
\label{}
\end{align}
The relation $ \le_{\rm icx} $ 
is a partial order on the class of random variables, 
so that it is called {\em increasing convex order} 
(see Shaked--Shanthikumar \cite{MR2265633}). 
It is known (Kellerer \cite{MR0356250}) that 
a process $ (S_t)_{t \ge 0} $ is increasing in increasing convex order 
if and only if 
$ (S_t)_{t \ge 0} $ has the same one-dimensional marginals with a submartingale. 
Interested readers are referred to 
Rothschild--Stiglitz (\cite{MR0503565},\cite{MR0503567}), 
Baker--Yor \cite{MR2519530}, and also Hirsch--Yor \cite{MR2571849}. 
\end{Rem}

\begin{Rem}
Profeta--Roynette--Yor \cite{MR2582990} 
proved that a Bessel process admits pseudo-inverse if and only if 
the dimension is greater than one, 
and investigated several remarkable properties of its pseudo-inverse. 
See also Yen--Yor \cite{YY} for another related study of Bessel process. 
\end{Rem}

This paper is organized as follows. 
In Section \ref{sec: dis}, we discuss d-inverses of several classes of processes 
and study scaling limit theorems of d-inverses. 
In Section \ref{sec: sca}, we study the inverse problem 
of scaling limits of d-inverses.

\section{Discussions on d-increasing processes} \label{sec: dis}

For two random variables $ X $ and $ Y $, 
we write $ X \dist Y $ if $ P(X \le x) = P(Y \le x) $ for all $ x \in \bR $. 
For a family of random variables $ (X^{(a)})_{a \in I} $ 
indexed by an interval $ I $ of $ \bR $, 
we write $ X^{(a)} \cdist X $ as $ a \to b \in I $ for a random variable $ X $ 
if $ P(X^{(a)} \le x) \to P(X \le x) $ as $ a \to b $ 
for all $ x \in \bR $ such that $ P(X=x)=0 $. 

\subsection{Transformations by increasing functions}

For an increasing function 
$ f:I \to[-\infty ,\infty ] $ 
defined on an subinterval $ I $ on $ \bR $, 
we denote its left-continuous inverse by $ f^{-1}:\bR \to [-\infty ,\infty ] $, 
i.e.: 
\begin{align}
f^{-1}(y) 
=& \inf \{ x \in I : f(x) \ge y \} 
\label{} \\
=& \sup \{ x \in I : f(x) < y \} , 
\label{}
\end{align}
where we adopt the usual convention that 
$ \inf \emptyset = \sup I $ and $ \sup \emptyset = \inf I $. 
By definition, we see that 
\begin{align}
f(x) \ge y \ \text{implies} \ x \ge f^{-1}(y) , 
\label{} \\
f(x) < y \ \text{implies} \ x \le f^{-1}(y) . 
\label{}
\end{align}

As a general remark, we give the following theorem. 

\begin{Thm} \label{thm1}
Let $ X=(X_t)_{t \ge 0} $ be a stochastic process 
such that $ X_t \in [x_0,\infty ) $ almost surely for all $ t \ge 0 $. 
Let $ f:[x_0,\infty ) \to \bR $ and $ g:[0,\infty ) \to [0,\infty ) $ 
be continuous increasing functions. 
Suppose that $ X $ admits a d-inverse $ (Y_x)_{x \ge x_0} $. 
Then $ \hat{X} = (\hat{X}_t)_{t \ge 0} $ defined by 
\begin{align}
\hat{X}_t = f \rbra{ X_{g(t)} } 
, \quad t \ge 0 
\label{}
\end{align}
admits a d-inverse $ \rbra{ g^{-1} \rbra{ Y_{f^{-1}(y)} } }_{y \ge f(x_0)} $. 
\end{Thm}

\Proof{
Since $ f $ is continuous and increasing, we see that 
$ f(f^{-1}(y))=y $, and hence that 
$ f(x) \ge y $ if and only if $ x \ge f^{-1}(y) $. 
This proves that 
\begin{align}
P(f(X_{g(t)}) \ge y) 
=& P(X_{g(t)} \ge f^{-1}(y)) 
\label{} \\
=& P(Y_{f^{-1}(y)} \le g(t)) 
\label{} \\
=& P \rbra{ g^{-1}(Y_{f^{-1}(y)}) \le t } . 
\label{}
\end{align}
The proof is complete. 
}

\subsection{Brownian motion with functional drift}

\begin{Thm} \label{thm2}
Let $ \rho:[0,\infty ) \to \bR $ be a right-continuous function. 
Then the process $ B^{(\rho)}_t = B_t + \rho(t) $ is d-increasing on $ [0,\infty ) $ 
if and only if the following condition is satisfied: 
\begin{align}
\text{\bf (A)} \quad 
\frac{\rho(t)}{\sqrt{t}} \ \text{is increasing in $ t>0 $}. 
\label{}
\end{align}
In this case, the d-inverse $ (Y^{(\rho)}_x)_{x \ge 0} $ is given by 
\begin{align}
Y^{(\rho)}_x \dist \eta_x^{-1}(B_1) 
\quad \text{for all $ x \ge 0 $}, 
\label{}
\end{align}
where $ \eta:(0,\infty ) \to \bR $ 
is the increasing function defined by 
\begin{align}
\eta_x(t) = \frac{ \rho(t)-x }{\sqrt{t}} 
, \quad t > 0 . 
\label{eq: Fx}
\end{align}
\end{Thm}

\Proof{
Since $ B_t \dist - \sqrt{t} B_1 $, we have 
\begin{align}
P \rbra{ B^{(\rho)}_t \ge x } 
= P \rbra{ B_1 \le \eta_x(t) } , 
\label{}
\end{align}
where $ \eta_x $ is defined as \eqref{eq: Fx}. 
Now $ B^{(\rho)} $ is d-increasing if and only if 
$ \eta_x(t) $ is increasing in $ t>0 $ for all $ x \ge 0 $, 
which is equivalent to the condition {\bf (A)}. 
}

In the remainder of this section, 
we discuss several particular classes of Browinan motion with functional drifts.

\subsection{Brownian motion with explosion}

Using $ B_t \dist \sqrt{t} B_1 $, we obtain the following: 
The Brownian motion without drift, $ B=B^{(0)} $, 
admits a d-inverse $ Y^{(0)}=(Y_x^{(0)})_{x \ge 0} $. 
In fact, we have 
\begin{align}
Y^{(0)}_x \dist 
\rbra{\frac{x}{B_1}}^2 1_{\{ B_1>0 \}} 
+ \infty 1_{\{ B_1 \le 0 \}} 
, \quad x \ge 0 . 
\label{}
\end{align}

For a constant $ t_0 \in (0,\infty ) $, 
the process $ X=(X_t)_{t \ge 0} $ 
taking values in $ (-\infty ,\infty ] $ defined by 
\begin{align}
X_t = B_t + \infty 1_{\{ t \ge t_0 \}} 
, \quad t \ge 0 
\label{}
\end{align}
is called {\em Brownian motion with explosion in finite time}. 
It admits a d-inverse $ Y=(Y_x)_{x \ge 0} $ given by 
\begin{align}
Y_x \dist \min \cbra{ Y^{(0)}_x , t_0 } 
, \quad x \ge 0 . 
\label{}
\end{align}

\begin{Thm}
Let $ \rho:[0,\infty ) \to (0,\infty ) $ 
be a right-continuous function satisfying the condition {\bf (A)}. 
Let $ \phi_1,\phi_2:[0,\infty ) \to [0,\infty ) $ be two functions. 
Suppose that there exist constants $ t_0 \in (0,\infty ] $ and $ p \in [0,\infty) $ such that 
\begin{align}
\text{\bf (B)} \quad 
\begin{cases}
\displaystyle \Bigg. 
\frac{\phi_1(\lambda) \rho(\lambda t)}{\sqrt{\lambda t}} 
\tend{}{\lambda \to 0+} 
\begin{cases}
0 & \text{if $ 0<t<t_0 $}, \\
\infty & \text{if $ t>t_0 $}, 
\end{cases}
\\
\displaystyle \Bigg. 
\frac{\phi_2(\lambda)}{\sqrt{\lambda }} 
\tend{}{\lambda \to 0+} p . 
\end{cases}
\label{}
\end{align}
Then, for any $ x \ge 0 $, it holds that 
\begin{align}
\frac{1}{\lambda} Y^{(\phi_1(\lambda) \rho)}_{\phi_2(\lambda) x} 
\cdist 
\min \cbra{ Y^{(0)}_{px} , t_0 } 
\quad \text{as $ \lambda \to 0+ $}. 
\label{eq: sc lim2}
\end{align}
In particular, for any $ \lambda > 0 $, it holds that 
\begin{align}
\frac{1}{\lambda} \min \cbra{ Y^{(0)}_{\sqrt{\lambda} x} , t_0 } 
\dist 
\min \cbra{ Y^{(0)}_x , t_0 } . 
\label{eq: sc inv2}
\end{align}
\end{Thm}

\Proof{
Since $ B_{\lambda t} \dist \sqrt{\lambda} B_t $, we have 
\begin{align}
P \rbra{ \frac{1}{\lambda} 
Y^{(\phi_1(\lambda) \rho)}_{\phi_2(\lambda) x} \le t } 
=& 
P \rbra{ B_{\lambda t} + \phi_1(\lambda) \rho(\lambda t) \ge \phi_2(\lambda) x } 
\label{} \\
=& 
P \rbra{ B_t + \frac{\phi_1(\lambda) \rho(\lambda t)}{\sqrt{\lambda}} 
\ge \frac{\phi_2(\lambda)}{\sqrt{\lambda}} x } . 
\label{}
\end{align}
The last quantity converges as $ \lambda \to 0+ $ 
to $ P(B_t \ge px) $ if $ t<t_0 $ and to $ 1 $ if $ t > t_0 $. 
Since we have 
\begin{align}
P \rbra{ \min \cbra{ Y^{(0)}_{px} , t_0 } \le t } 
= 
\begin{cases}
P(B_t \ge px) & \text{if $ t<t_0 $}, \\
1 & \text{if $ t \ge t_0 $}, 
\end{cases}
\label{}
\end{align}
we obtain \eqref{eq: sc lim2}. 
The scale invariance property \eqref{eq: sc inv2} is obvious. 
The proof is now complete. 
}

\subsection{Brownian motion with constant drift}

By Theorem \ref{thm2}, we see that 
the Brownian motion with constant drift $ B^{(c \cdot)} = (B_t + ct)_{t \ge 0} $ 
admits a d-inverse $ Y^{(c \cdot)}=(Y_x^{(c \cdot)})_{x \ge 0} $ 
if and only if $ c \in [0,\infty) $. 
If $ c \in (0,\infty) $, i.e., except for the Brownian case, 
we obtain, for $ x \ge 0 $, 
\begin{align}
Y_x^{(c \cdot)} \dist 
\rbra{ \frac{B_1 + \sqrt{B_1^2 + 4cx}}{2c} }^2 . 
\label{}
\end{align}
We remark that, for any $ x \ge 0 $, 
\begin{align}
Y_x^{(c \cdot)} \cdist Y_x^{(0)} 
\quad \text{as $ c \to 0+ $}. 
\label{}
\end{align}
We also remark the following: 
Using $ B_t \dist - t B_{1/t} $, we can easily see that 
\begin{align}
Y^{(c \cdot)}_x \dist \frac{1}{Y^{(x \cdot)}_c} 
\quad \text{for all $ c \ge 0 $ and $ x \ge 0 $}. 
\label{}
\end{align}

Scaling property of Brownian motion with constant drifts will be discussed 
in the next section in a more general setting. 

The geometric Brownian motion $ S=S^{(\sigma,\mu)} $ 
with constant volatility $ \sigma > 0 $ and drift $ \mu \in \bR $ 
given as \eqref{eq: gbm} 
may be represented as $ S^{(\sigma,\mu)}_t = f(B^{((\tilde{\mu}/\sigma) t)}_t) $ where 
$ f(x) = s_0 \exp(\sigma x) $. 
Hence we may apply Theorem \ref{thm1} and obtain the following: 
$ S^{(\sigma,\mu)} $ admits a d-inverse $ (T^{(\sigma,\mu)}_s)_{s \ge s_0} $ 
if and only if $ \tilde{\mu} = \mu - \sigma^2/2 \ge 0 $. 
In this case, we have 
\begin{align}
T^{(\sigma,\mu)}_s 
\dist 
Y^{((\tilde{\mu}/\sigma) \cdot)}_{f^{-1}(s)} 
\quad \text{for all $ s \ge s_0 $}. 
\label{}
\end{align}

\subsection{Brownian motion with power drift}

For $ \alpha \in [0,\infty) $ and $ c \in [0,\infty ) $, we define 
\begin{align}
R^{(c,\alpha )}_t = B_t + ct^{\alpha } 
, \quad t \ge 0 
\label{}
\end{align}
and we call $ R^{(c,\alpha )} = (R^{(c,\alpha )}_t)_{t \ge 0} $ 
a {\em Brownian motion with power drift}. 
By Theorem \ref{thm2}, we see that 
$ R^{(c,\alpha )} $ admits a d-inverse $ (Z^{(c,\alpha )}_x)_{x \ge 0} $ 
if and only if $ \alpha \ge 1/2 $. 

The following theorem tells us that 
the class of the d-inverses of Brownian motion with power drifts 
appear as scaling limits, 
and consequently, satisfy scale invariance property. 

\begin{Thm} \label{thm3}
Let $ \rho:[0,\infty ) \to (0,\infty ) $ 
be a right-continuous function satisfying the condition {\bf (A)}. 
Let $ \phi_1,\phi_2:[0,\infty ) \to [0,\infty ) $ be two functions. 
Suppose there exist $ \alpha \ge 1/2 $, $ c\in (0,\infty) $ and $ p \in [0,\infty) $ such that 
\begin{align}
\text{\bf (RV)} \quad 
\begin{cases}
\displaystyle \Bigg. 
\frac{\rho(\lambda t)}{\rho(\lambda)} \tend{}{\lambda \to 0+} t^{\alpha } , \\
\displaystyle \Bigg. 
\frac{\rho(\lambda)}{\sqrt{\lambda}} \phi_1(\lambda) \tend{}{\lambda \to 0+} c , \\
\displaystyle \Bigg. 
\frac{1}{\sqrt{\lambda}} \phi_2(\lambda) \tend{}{\lambda \to 0+} p . 
\end{cases}
\label{eq: rv}
\end{align}
Then, for any $ x \ge 0 $, it holds that 
\begin{align}
\frac{1}{\lambda} Y^{(\phi_1(\lambda) \rho)}_{\phi_2(\lambda) x} 
\cdist 
Z^{(c,\alpha )}_{px} 
\quad \text{as $ \lambda \to 0+ $}. 
\label{eq: sc lim}
\end{align}
In particular, for any $ \lambda > 0 $, it holds that 
\begin{align}
\frac{1}{\lambda} Z^{\rbra{ c \lambda^{(1/2)-\alpha },\alpha }}_{\sqrt{\lambda} x} 
\dist 
Z^{(c,\alpha )}_x . 
\label{eq: sc inv}
\end{align}
\end{Thm}

\begin{Rem}
The condition {\bf (RV)} asserts that 
the functions $ \rho $, $ \phi_1 $ and $ \phi_2 $ (if $ p\in (0,\infty) $) 
are regularly varying at $ 0+ $ of index 
$ \alpha $, $ (1/2)-\alpha $, and $ 1/2 $, respectively. 
\end{Rem}

\Proof[Proof of Theorem \ref{thm3}]{
Since $ B_{\lambda t} \dist \sqrt{\lambda} B_t $, we have 
\begin{align}
P \rbra{ \frac{1}{\lambda} 
Y^{(\phi_1(\lambda) \rho)}_{\phi_2(\lambda) x} \le t } 
=& 
P \rbra{ B_t + \frac{\rho(\lambda)}{\sqrt{\lambda}} \phi_1(\lambda) 
\cdot \frac{\rho(\lambda t)}{\rho(\lambda)} \ge \frac{\phi_2(\lambda)}{\sqrt{\lambda}} x } 
\label{} \\
\tend{}{\lambda \to 0+}& 
P(B_t + c t^{\alpha } \ge px) 
\label{} \\
=& P(Z^{(c,\alpha )}_{px} \le t) . 
\label{}
\end{align}
Now we have obtained \eqref{eq: sc lim}. 
The scale invariance property \eqref{eq: sc inv} is obvious. 
The proof is complete. 
}

\section{Scaling limits for the class of d-inverses} \label{sec: sca}

In what follows, by {\em measurable} we mean Lebesgue measurable. 

\begin{Thm} \label{thm4}
Let $ \rho:[0,\infty ) \to [0,\infty ) $ 
be a right-continuous function satisfying the condition {\bf (A)}. 
Suppose that, 
for some measurable functions $ \phi_1,\phi_2:(0,\infty ) \to (0,\infty ) $ 
and for some family $ Z=(Z_x)_{x \ge 0} $ of $ [0,\infty ] $-valued random variables, 
it holds that 
\begin{align}
\frac{1}{\lambda} Y^{\rbra{\phi_1(\lambda) \rho}}_{\phi_2(\lambda) x} 
\tend{\rm d}{\lambda \to 0+} 
Z_x 
\quad \text{for all $ x \ge 0 $}. 
\label{eq: ass}
\end{align}
Suppose, moreover, that there exists a constant $ p \in [0,\infty) $ such that 
\begin{align}
\frac{\phi_2(\lambda)}{\sqrt{\lambda}} \tend{}{\lambda \to 0+} p . 
\label{eq: ass2}
\end{align}
Then either one of the following four assertions holds: 
\begin{enumerate}
\item 
$ \phi_1(\lambda) \rho(\lambda t)/\sqrt{\lambda} \tend{}{\lambda \to 0+} 0 $ for all $ t>0 $. 
In this case, 
\begin{align}
Z_x \dist Y^{(0)}_{px} 
\quad \text{for all $ x \ge 0 $}. 
\label{eq: 0 infty limit3-}
\end{align}
\item 
The condition {\bf (B)} holds 
for some $ t_0 \in (0,\infty ) $. 
In this case, 
\begin{align}
Z_x \dist \min \cbra{ Y^{(0)}_{px} , t_0 } 
\quad \text{for all $ x \ge 0 $}. 
\label{eq: 0 infty limit3}
\end{align}
\item 
The condition {\bf (RV)} holds 
for some $ \alpha \ge 1/2 $ and $ c\in (0,\infty) $. 
In this case, 
\begin{align}
Z_x \dist Z^{(c,\alpha )}_{px} 
\quad \text{for all $ x \ge 0 $}. 
\label{}
\end{align}
\item 
(Degenerate case.) 
$ P(Z_x=0)=1 $ for all $ x\in (0,\infty) $. 
\end{enumerate}
\end{Thm}

\Proof{
Let $ x \ge 0 $. 
Denote $ F_x(t) = P(Z_x \le t) $ for $ t \ge 0 $ 
and denote by $ C(F_x) $ the set of continuity point of $ F_x $. 
We note that 
\begin{align}
P \rbra{ \frac{1}{\lambda} Y^{\rbra{\phi_1(\lambda) \rho}}_{\phi_2(\lambda) x} \le t } 
=& P \rbra{ B_{\lambda t} + \phi_1(\lambda) \rho(\lambda t) \ge \phi_2(\lambda) x } 
\label{} \\
=& P \rbra{ B_1 + \phi_1(\lambda) \frac{\rho(\lambda t)}{\sqrt{\lambda t}} 
\ge \frac{\phi_2(\lambda)}{\sqrt{\lambda t}} x } . 
\label{}
\end{align}
By the assumption \eqref{eq: ass}, we see that 
\begin{align}
\begin{split}
P \rbra{ B_1 + \phi_1(\lambda) \frac{\rho(\lambda t)}{\sqrt{\lambda t}} 
- \frac{\phi_2(\lambda)}{\sqrt{\lambda t}} x \in [0,\infty) } 
\tend{}{\lambda \to 0+} 
P(Z_x \le t) 
\\
\quad \text{for all $ t \in C(F_x) \cap (0,\infty ) $}. 
\end{split}
\label{eq: psi1 rho limit-}
\end{align}
Hence there exists a function $ g_x: C(F_x) \cap (0,\infty ) \to [-\infty ,\infty ] $ such that 
\begin{align}
\phi_1(\lambda) \frac{\rho(\lambda t)}{\sqrt{\lambda t}} 
- \frac{\phi_2(\lambda)}{\sqrt{\lambda t}} x 
\tend{}{\lambda \to 0+} 
g_x(t) 
\quad \text{for all $ t \in C(F_x) \cap (0,\infty ) $}. 
\label{eq: psi1 rho limit}
\end{align}
Since $ \rho $ satisfies the condition {\bf (A)} 
and since $ C(F_x) $ is dense in $ \bR $, 
we see that $ g_x $ is increasing, 
and hence we may extend $ g_x $ on $ [0,\infty ) $ so that it is right-continuous. 
Now we obtain, for any $ x \ge 0 $, 
\begin{align}
Z_x \dist g_x^{-1}(B_1) . 
\label{eq: Zx dist}
\end{align}

Let us write $ g $ simply for $ g_0 $. 
Noting that $ g $ is an increasing function taking values in $ [0,\infty ] $. 
we divide into the following four distinct cases. 

(i) {\em The case where $ g(t)=0 $ for all $ t>0 $.} 

Let $ x \ge 0 $ be fixed. By the assumption \eqref{eq: ass2} 
and by \eqref{eq: psi1 rho limit}, we obtain 
\begin{align}
g_x(t) = - px/\sqrt{t} , \quad t>0 . 
\label{}
\end{align}
From this and \eqref{eq: Zx dist}, we obtain 
\begin{align}
P(Z_x \le t) = P(Y^{(0)}_{px} \le t) 
, \quad t>0 . 
\label{}
\end{align}
This proves \eqref{eq: 0 infty limit3-}. 
The proof of Claim (i) is now complete. 

(ii) {\em The case where there exist a point $ t_0 \in (0,\infty ) $ such that 
\begin{align}
g(t) 
\begin{cases}
= 0 & \text{if $ 0<t<t_0 $}, \\
= \infty & \text{if $ t>t_0 $}. 
\end{cases}
\label{}
\end{align}}

Let $ x \ge 0 $. By the assumption \eqref{eq: ass2} 
and by \eqref{eq: psi1 rho limit}, we obtain 
\begin{align}
g_x(t) = 
\begin{cases}
- px/\sqrt{t} & \text{if $ 0<t<t_0 $}, \\
\infty & \text{if $ t>t_0 $}. 
\end{cases}
\label{}
\end{align}
From this and \eqref{eq: Zx dist}, we obtain 
\begin{align}
P(Z_x \le t) = 
\begin{cases}
P(Y^{(0)}_{px} \le t) & \text{if $ 0 \le t < t_0 $}, \\
1 & \text{if $ t \ge t_0 $}. 
\end{cases}
\label{}
\end{align}
This proves \eqref{eq: 0 infty limit3}. 
The proof of Claim (ii) is now complete. 

(iii) {\em The case where 
there are two points $ t_0 , t_1 \in C(F_0) \cap (0,\infty ) $ with $ t_0<t_1 $ 
such that $ 0 < g(t_0) \le g(t_1) < \infty $.} 

Since $ g $ is increasing, we see that 
\begin{align}
0<g(t)<\infty \quad \text{for all $ t \in C(F_0) \cap [t_0,t_1] $}. 
\label{}
\end{align}
By \eqref{eq: psi1 rho limit}, we have, for any $ t \in C(F_0) \cap [t_0,t_1] $, 
\begin{align}
\frac{\rho(\lambda t)}{\rho(\lambda t_0)} 
= \frac{\phi_1(\lambda) \frac{\rho(\lambda t)}{\sqrt{\lambda t}} }
{\phi_1(\lambda) \frac{\rho(\lambda t_0)}{\sqrt{\lambda t_0}} } 
\cdot \frac{\sqrt{t}}{\sqrt{t_0}} 
\tend{}{\lambda \to 0+} 
\frac{g(t)}{g(t_0)} \cdot \frac{\sqrt{t}}{\sqrt{t_0}} \in (0,\infty ) . 
\label{eq: rho ratio lim}
\end{align}
Since $ C(F_0) \cap [t_0,t_1] $ has positive Lebesgue measure, 
we may apply Characterisation Theorem (\cite[Theorem 1.4.1]{MR898871}) to see that 
the convergence \eqref{eq: rho ratio lim} 
and consequently \eqref{eq: psi1 rho limit} are still valid for all $ t \in (0,\infty ) $, 
and that 
\begin{align}
\frac{g(t)}{g(t_0)} \cdot \frac{\sqrt{t}}{\sqrt{t_0}} = t^{\alpha } 
, \quad t \in (0,\infty ) 
\label{}
\end{align}
for some $ \alpha \in \bR $. 
Since $ g $ is increasing, we have $ \alpha \ge 1/2 $. 
We obtain 
\begin{align}
g(t) = c t^{\alpha -1/2} 
, \quad t \in (0,\infty ) 
\label{}
\end{align}
for some $ c \in (0,\infty ) $. 
Hence, by \eqref{eq: rho ratio lim} and \eqref{eq: psi1 rho limit}, we obtain 
\begin{align}
\frac{\rho(\lambda t)}{\rho(\lambda)} \tend{}{\lambda \to 0+} t^{\alpha } 
\quad \text{and} \quad 
\frac{\rho(\lambda)}{\sqrt{\lambda}} \phi_1(\lambda) \tend{}{\lambda \to 0+} c . 
\label{}
\end{align}
Now we have seen that the condition {\bf (RV)} is satisfied. 
The proof of Claim (iii) is now completed by Theorem \ref{thm3}. 

(iv) {\em The case where $ g(t) = \infty $ for all $ t>0 $.} 

In this case, by the assumption \eqref{eq: ass2} 
and by \eqref{eq: psi1 rho limit}, we obtain 
$ g_x(t)=\infty $ for all $ t>0 $ and $ x \ge 0 $. 
By \eqref{eq: Zx dist}, we obtain $ P(Z_x = 0) = 1 $ for all $ x \ge 0 $. 
The proof of Claim (iv) is now complete. 
}

\end{document}